\documentclass[
final
 , nomarks
]{dmtcs-episciences}

\usepackage[utf8]{inputenc}
\usepackage{subfigure}
\usepackage[centertags,intlimits]{amsmath}

\newtheorem{theorem}{Theorem}[section]
\newtheorem{proposition}[theorem]%
{Proposition}
\newtheorem{definition}{Definition}[section]

\newtheorem{lemma}[theorem]%
{Lemma}
\newtheorem{corollary}[theorem]%
{Corollary}

\usepackage{hyperref}

\author{Peter Dankelmann\affiliationmark{1}
  \and Sonwabile Mafunda\affiliationmark{1}
  \and Sufiyan Mallu\affiliationmark{1}\thanks{Financial support by the South African National Research Foundation
is gratefully acknowledged.}}
\title[Proximity, remoteness and maximum degree in graphs]{Proximity, remoteness and maximum degree \\ in graphs}
\affiliation{
  University of Johannesburg, South Africa}
\received{2022-05-06}
\revised{2022-10-22}
\accepted{2022-10-24}
\begin{document}
\publicationdetails{24}{2022}{2}{10}{9432}
\maketitle
\begin{abstract}
The average distance of a vertex $v$ of a connected graph $G$ is the arithmetic mean of the distances 
from $v$ to all other vertices of $G$. The proximity $\pi(G)$ and the remoteness $\rho(G)$ of $G$ 
are the minimum and the maximum of the average distances of the vertices of $G$, respectively.

In this paper, we give upper bounds on the remoteness and proximity for graphs of given order,
minimum degree and maximum degree. Our bounds are sharp apart from an additive constant. 
\end{abstract}

\begin{center}
{\it To the memory of Aisha Patel}
\end{center}

\section{Introduction}
\label{sec:in}

Let $G$ be a connected graph of order $n$ with vertex set $V(G)$. The {\em average distance} 
$\overline{\sigma}(v)$ of a vertex $v$ of $G$ is defined as the arithmetic mean of the 
distances from $v$ to all other vertices of $G$, i.e., 
\[ \overline{\sigma}(v,G) = \frac{1}{n-1} \sum_{u \in V(G)} d(v,u),   \]
where $d(v,u)$ is the usual shortest path distance between vertices $v$ and $u$. 
The {\em proximity} and the {\em  remoteness} of $G$, denoted by $\pi(G)$ and $\rho(G)$, are the
smallest and the largest, respectively, average distance among the vertices of $G$. 
The average distance of a vertex $v$ is closely related to its {\em total distance} $\sigma(v,G)$, 
defined as the sum of the distances from $v$ to all other vertices. 
Clearly, $\sigma(v,G) = (n-1) \overline{\sigma}(v,G)$. 
Also the names {\em status} or {\em transmission} have been used in the
literature. The proximity of a graph is closely related to its {\em minimum status},
defined as the smallest total distance among the vertices of the graph. 

If a graph $G$ represents a network in which we want to place a facility which should be 
close, on average, to vertices of $G$, then the proximity of $G$ is an indicator for 
how good the best location in the network is, and the remoteness is an indicator for how
good the worst location is.

The proximity and also the remoteness of a connected graph on at least two vertices is at least $1$. 
Sharp upper bounds in terms of order alone were given by 
Zelinka \cite{Zel1968} and later, independently, by Aouchiche and Hansen \cite{AouHan2011}, who
introduced the names proximity and remoteness. They proved that for every connected graph 
of order $n$, 
\begin{equation}  \label{eq:remoteness-in-terms-of-order}
\rho(G) \leq \frac{n}{2}, 
\end{equation}
with equality if and only if $G$ is a path, and also that 
\begin{equation} \label{eq:proximity-in-terms-of-order}  
\pi(G) \leq  \left\{ \begin{array}{cc}
\frac{n+1}{4} & \textrm{if $n$ is odd,} \\
\frac{n+1}{4} + \frac{1}{4(n-1)} & \textrm{if $n$ is even,} 
\end{array} \right. 
\end{equation}
with equality if and only if $G$ is a path or a cycle.

There are several results in the literature on relations between proximity or remoteness 
and other distance measures. 
For the {\em diameter} (defined as the largest of the distances
between the vertices), these were first investigated by 
Aouchiche and Hansen \cite{AouHan2011}, who determined sharp upper bounds on the difference 
between diameter and proximity and on the difference between diameter and remoteness in terms
of order. Improved bounds, that take into account also the minimum degree were 
given in \cite{Dan2016} and \cite{DanMaf-manu}.

Also bounds involving the {\em radius} (defined as the smallest of the eccentricities of the 
vertices of $G$, where the eccentricity of a vertex $v$ is the distance from $v$ to a vertex 
farthest from $v$) have been explored. A sharp upper bound on the  
difference between radius and proximity of a graph of given order was given in \cite{AouHan2011}, 
and improved bounds that take into account also the minimum degree can be found in 
\cite{Dan2016} and \cite{DanMaf-manu}.
A conjecture in \cite{AouHan2011} on the maximum value of the difference between remoteness and 
radius for graphs of given order was proved independently by Wu and Zhang \cite{WuZha2014} and 
Hua, Chen and Das \cite{HuaCheDas2015}.

Ma, Wu and Zhang \cite{MaWuZha2012} showed that the difference between {\em average eccentricity} 
(defined as the arithmetic mean of the eccentricities of the vertices) 
and proximity is maximised by the path, thus proving a conjecture from \cite{AouHan2011}. 
Another conjecture in \cite{AouHan2011} which states that the difference between average
eccentricity and remoteness cannot exceed that of the cycle was shown by Sedlar \cite{Sed2013}
to hold for trees. 

The maximum value of the difference between remoteness and {\em average distance} (defined as 
the arithmetic 
mean of the distances between all vertices of the graph) and the difference
between average distance and proximity were determined by Wu and Zhang 
\cite{WuZha2014} and by Sedlar \cite{Sed2013}, respectively, thus confirming two
conjectures from \cite{AouHan2011}. The minimum value of the ratio of proximity to
average distance for graphs of given order was determined by Hua and Das \cite{HuaDas2014}. 

Proximity and remoteness have been studied for graphs from various classes. Among other results, 
Barefoot, Entringer and Sz\'{e}kely \cite{BarEntSze1997} determined the maximum value of 
the ratio $\rho(T)/\pi(T)$ for a tree $T$ of given order. They also determined lower bounds 
on the ratios of average distance to proximity and 
average distance to remoteness for trees of given order. 
The trees that maximise proximity among trees of given order and diameter were determined 
by Peng and Zhou \cite{PenZho2021}. The same authors also gave bounds on proximity in terms
of order and either number of end-vertices, number of vertices of odd degree and number
of vertices of degree $2$. 
Proximity of series-reduced trees, i.e., trees with no vertex of degree $2$, were 
studied by Cheng, Lin and Zhou \cite{CheLinZho2021}. 
For maximal planar graphs, bounds on remoteness and proximity were given by Czabarka, Dankelmann, 
Olsen and Sz\'{e}kely in \cite{CzaDanOlsSze2019} and \cite{CzaDanOlsSze2020}, respectively. 
The study of proximity and remoteness in digraphs was initiated by Ai, Gerke, Gutin and 
Mafunda \cite{AiGerGutMaf2021}.

This paper is concerned with bounds on proximity and remoteness that take into account
vertex degrees. The {\em degree} of a vertex $v$ is defined by ${\rm deg}_G(v) =|N_G(v)|$,
where $N_G(v)$ is the {\em neighbourhood} of $v$ i.e., the set of vertices adjacent
to $v$. By $\delta(G)$ and $\Delta(G)$ we denote the {\em minimum degree} and the 
{\em maximum degree} of $G$, i.e., the smallest and the largest degree of a vertex
in $G$. 
 
The upper bounds in \eqref{eq:remoteness-in-terms-of-order} and 
\eqref{eq:proximity-in-terms-of-order} are attained by paths, which have minimum
degree $1$. For graphs of  larger minimum degree $\delta$, the following improved bounds 
were given in \cite{Dan2015}.
\begin{equation} \label{eq:remoteness-vs-minimum-degree}
\rho(G) \leq \frac{3n}{2(\delta+1)} +\frac{7}{2}, 
\end{equation}
\begin{equation}  \label{eq:proximity-vs-minimum-degree}
\pi(G) \leq \frac{3n}{4(\delta+1)} +3, 
\end{equation}
and further improvements  for graphs not containing a $3$-cycle or a $4$-cycle as a subgraph were given in \cite{DanJonMaf2021}. 

The graphs constructed in \cite{Dan2015} to show that \eqref{eq:remoteness-vs-minimum-degree} 
and \eqref{eq:proximity-vs-minimum-degree} are sharp apart from an additive constant  
are close to regular. Hence it is natural to ask if one can find improved bounds 
for graphs containing a vertex of large degree. That this is indeed the case for trees was
shown by Tsai, Shang and Zhang \cite{LinTsaShaZha2012} and Rissner and Burkhard \cite{RisBur2014},
who determined the trees of given order and maximum degree that maximise proximity and remoteness.
In this paper we answer the above question in the affirmative and improve the 
bounds \eqref{eq:remoteness-vs-minimum-degree} and 
\eqref{eq:proximity-vs-minimum-degree} for graphs of given maximum degree. We prove the 
following bounds, which are sharp apart from an additive constant: 
\[ \pi(G) \leq \left\{ \begin{array}{cc} 
\frac{3(n-\Delta)^2}{2(n-1)(\delta+1)} + \frac{13}{2} & \textrm{if $\Delta > \frac{n}{2}-1$,} \\[1mm] 
 \frac{3n^2 - 6\Delta^2}{4(n-1)(\delta+1)} + \frac{35}{4} 
   & \textrm{if $\Delta  \leq \frac{n}{2}-1$}, \end{array}  \right. \]
and 
\[ \rho(G) \leq  \frac{3(n^2 -\Delta^2)}{2(n-1)(\delta+1)} + 7, \]
where $\Delta$ denotes the maximum degree of $G$.  
Our bounds show a certain analogy to results in \cite{AloDan2021}, where it was shown that 
bounds on the average distance of graphs in terms of order and minimum degree can be improved
significantly for graphs with large maximum degree.

This paper is organised as follows. In 
Section \ref{section:Proximity and remoteness of vertex-weighted graphs} we consider 
graphs with a weight function on the vertex set. We define the weighted distance of a vertex 
and prove bounds on the weighted distance. These bounds are used in 
Section \ref{section:Proximity and remoteness of graphs of given minimum degree and maximum degree}
to prove our main results, bounds on proximity and remoteness in terms of order, minimum degree
and maximum degree. 
Graphs that show that these bounds are sharp  apart from an additive constant are constructed 
in Section \ref{section:A sharpness example}.

\section{Proximity and remoteness of weighted graphs}
\label{section:Proximity and remoteness of vertex-weighted graphs}

In this section we consider graphs with a weight function on the vertex set. 
Taking an approach similar to
that in \cite{Dan2015}, we first define the weighted distance of a vertex $v$. 

\begin{definition}
Let $G$ be a connected graph and $c: V(G) \longrightarrow \reals^{\geq 0}$ be a 
nonnegative weight function on the vertices of $G$. Let $v$ be a vertex of $G$. Then the 
weighted distance of $v$ with respect to $c$ is defined as 
\[ \sigma_c(v) = \sum_{w \in V(G) -\{v\}} c(w) d(v,w). \]
A vertex whose weighted distance is minimum among all vertices of $G$ is a $c$-median vertex,
and the $c$-median of $G$ is the set of all $c$-median vertices of $G$.
\end{definition}

If $c(v)=1$ for every vertex $v$ of $G$, then the $c$-median is exactly the median of $G$. 
It is well-known (see for example \cite{Zel1968}) that in a tree $T$  of order $n$, 
the median vertices are exactly
the vertices whose branch weight is not more than $\frac{n}{2}$, where the branch weight 
of a vertex $v$ in $T$ is defined as the maximum order among all components of $T-v$. 
A more general result for weighted trees was given by Kariv and Hakimi \cite{KarHak1979}. 

\begin{definition}
Let $T$ be a tree and $c: V(T) \longrightarrow \reals^{\geq 0}$ be a 
nonnegative weight function on the vertices of $T$. Let $v$ be a vertex of $T$.
Then the $c$-branch weight ${\rm bw}_c(v,T)$ of $v$ is the largest 
weight of a component of $T-v$, where the weight of a component is the sum
of the weights of its vertices.
\end{definition}

If $G$ is a graph with a weight function $c$ on the vertex set, then for a set 
$A$ of vertices of $G$ we write $c(A)$ for $\sum_{v\in A} c(v)$. If $H$ is a 
subgraph of $G$, then we write $c(H)$ for $c(V(H))$.

\begin{proposition}[\cite{KarHak1979}]  \label{prop:median-branchweight-weighted} 
Let $T$ be a tree and $c: V(T) \longrightarrow \reals^{\geq 0}$ be a 
nonnegative weight function on the vertices of $T$. Let $v$ be a vertex of $T$.
Then $v$ is a $c$-median vertex of $T$ if and only if ${\rm bw}_c(v,T) \leq \frac{c(T)}{2}$. 
\end{proposition}

A key result on weighted graphs in \cite{Dan2015} gives a 
bound on the weighted distance of a $c$-median vertex of a weighted graph with given 
total weight but no restriction on the order, in which the weight of every vertex is
at least a prescribed value $k$. In the following lemma, which can be viewed as an extension 
of the result in \cite{Dan2015}, the graph satisfies the additional condition that the weight
of at least one vertex is not less than a prescribed large value $L$. Its proof is 
significantly more involved than that of the corresponding result in \cite{Dan2015}.

\begin{lemma}  \label{la:proximity-of-weighted}
Let $G$ be a connected graph and $k,L \in \reals$ with $0 < k < L$. 
Let $c: V(G) \longrightarrow \reals^{\geq 0}$ be
a weight function with total weight $N$. Assume that $c(u) \geq k$ for 
every vertex $u \in V(G)$, and that $G$ contains a vertex of weight at least $L$. 
Assume further that $N-L$ is an integer multiple of $k$. If $v$ is a $c$-median vertex of $G$, then
\begin{equation} \label{eq:statement-of-main-lemma}
\sigma_c(v,G) \leq \left\{ \begin{array}{cc}
   \frac{(N-L)(N-L+k)}{2k} & \textrm{if $L > \frac{N}{2}$,} \\[1mm]
   \frac{N^2- 2L^2}{4k} + \frac{N+L}{2}  & \textrm{if $L \leq \frac{N}{2}$.}
        \end{array} \right. 
\end{equation}        
\end{lemma}

{\bf Proof:}
Assume that $N$, $k$ and $L$ are given. 
Let $T$ be a spanning tree of $G$ that preserves the distances from a $c$-median vertex
$v$. Clearly, $v$ is a $c$-median vertex of $T$ and $\sigma_c(v,T)=\sigma_c(v,G)$, so it suffices 
to prove the lemma for $T$. We may assume that 
$T$ and $c$ are such that the weighted distance of a $c$-median vertex of $T$ is maximum among 
all trees and weight functions satisfying the hypothesis of the lemma. 

Our proof strategy is as follows. We prove a sequence of claims from which it follows 
that $T$ is a path in which one end-vertex has weight $L$, and all other vertices have weight
$k$. Evaluating the weighted distance of a $c$-median vertex of this path then yields the inequality \eqref{eq:statement-of-main-lemma}. 

We assume that $v$ is a $c$-median vertex of $T$ that is also an internal vertex, if possible. Let
$u_1, u_2,\ldots, u_d$ be the neighbours of $v$. For $i \in \{1,2, . . . ,d\}$ denote 
the component of $T-v$  containing $u_i$ by $T_i$, and its total weight by $C_i$.   
We may assume that $C_1 \geq C_2 \geq \ldots \geq C_d$. It follows from 
Proposition \ref{prop:median-branchweight-weighted} that 
$C_1 = {\rm bw}_c(v,T)  \leq  \frac{N}{2}$. \\[1mm]
{\sc Claim 1:} $C_i + C_j > \frac{N}{2}$ for all distinct $i,j \in \{1,2,\ldots,d\}$. \\
Suppose to the contrary that $C_i + C_j \leq \frac{N}{2}$
for some distinct $i,j \in \{1,2, \ldots, d\}$.
We consider the tree $H =T-vu_i + u_iu_j$ with the same weight function $c$. 
Clearly, the weights of the vertices of $H$ satisfy the hypothesis of the lemma. Furthermore, 
$v$ is also a $c$-median vertex of $H$ since the new branch has weight $C_i + C_j \leq \frac{N}{2}$, 
which implies that the $c$-branch weight of $v$ in $H$ is still at most $\frac{N}{2}$. 
Since the distance between the vertices of $T_i$ and $v$ has increased by $1$, we have
\[ \sigma_c(v,H) = \sigma_c(v,T) + \sum_{x \in V(T_i)} c(x) = \sigma_c(v,T) + C_i > \sigma_c(v,T), \]
contradicting our choice of $T$ and $v$. This proves Claim 1. \\[1mm]
{\sc Claim 2:} $T-v$ has at most $3$ components, i.e., $d \leq 3$. \\[1mm]
Suppose to the contrary that $T-v$ has at least $4$ components. Since by Claim 1 the total 
weight of any two of these is greater than $\frac{N}{2}$, the total weight of these four 
components would exceed $N$. This contradiction proves Claim 2. \\[1mm]
{\sc Claim 3:} $T_i$ is a path and $v$ is adjacent to one of the end-vertices of $T_i$ 
for each $i \in \{1,2,\ldots,d\}$. \\[1mm]
It suffices to prove that no vertex $w$ of $T_i$ has two neighbours that are farther by $1$ from
$v$ than $w$. Suppose to the contrary that some vertex $w$ of $T_i$ has two neighbours, say, 
$x_1$ and $x_2$ with $d(v,x_1) = d(v,x_2) = d(v,w)+1$. 
Consider the tree $H = T- wx_1 +x_1x_2$ with the same weight function $c$. The distance between 
$x_1$ and $v$ has increased by $1$, and no distance between $v$ and another vertex has decreased. 
Furthermore, $v$ is a $c$-median vertex of $H$ since its $c$-branch weight has not changed 
and is thus still at most $\frac{N}{2}$. Hence $\sigma_c(v,H) > \sigma_c(v,T)$. This 
contradiction to our choice of $T$ proves Claim 3. \\[-3mm]

For the remainder of the proof we use the following notation. 
For $i \in \{1,2,\ldots,d\}$ denote the end-vertex of $T_i$ that is farthest from $v$ by $z_i$. 
Let $y$ be a vertex of maximum weight in $T$, so $c(y) \geq  L > k$. If $v$ is among the vertices 
of maximum weight, then choose $y = v$.\\[1mm]
{\sc Claim 4:} Let $w$ be an internal vertex of $T$. Then $c(w) = k$, unless $w = v = y$ and
$c(w) = L$. \\
Assume that $T$ contains an internal vertex $w$ with $c(w) > k$. It suffices to show that
$w = v$, $w = y$ and $c(y)=L$. First suppose that $w \neq v$. Then $w \in V(T_i)$ for some 
$i \in \{1,2, \ldots, d\}$.
We obtain a weight function $c'$ from $c$ by reducing the weight of vertex $w$ to $k$ and
adding the difference $c(w)-k$ to the weight of $z_i$. Then $v$ is also a $c'$-median 
vertex since its $c'$-branch weight equals its $c$-branch weight, which is 
at most $\frac{N}{2}$. 
The conditions on the weight function are satisfied by $c'$. Indeed, every vertex
has weight at least $k$, and either $w \neq y$, in which case $c'(y) \geq L$, or $w=y$, in which
case $c'(z_i) = c(z_i)+c(y)-k \geq  k+L-k =L$, so there exists a vertex of weight at least 
$L$. Clearly, we have
\[ \sigma_{c'}(v,T) = \sigma_c(v,T) + d_T(w,z_i) \big(c(w)-k\big) > \sigma_c(v,T), \]
contradicting our choice of $T$ and $c$. This proves that $w=v$. 

Now suppose that $w \neq y$ or that $w = y$ and $c(w) > L$. Since $v = w$ and $w$ is
an internal vertex, $T-v$ has at least two components. Note that 
$C_1 + C_2 \leq N-c(v) < N-k$, thus we get that $C_2 \leq \frac{N}{2} - \frac{k}{2}$. 
We now obtain the weight function $c'$ by reducing
the weight of $w$ by $\min\{ \frac{k}{2}, c(w)-k\}$  if $w \neq y$, and by 
$\min\{ \frac{k}{2}, c(w)-L\}$ if $w = y$ and
increasing the weight of $z_2$ by the same amount. Then $v$ is also a $c'$-median vertex, $c'$
satisfies the hypothesis of the lemma, but $\sigma_{c'}(v,T) > \sigma_c(v,T)$, a 
contradiction to the choice of $T$ and $c$. We conclude that $w = y$ and $c(y)=L$, so 
Claim 4 follows. \\[1mm]
{\sc Claim 5:} The lemma holds if $L > \frac{N}{2}$. \\
Assume that $L > \frac{N}{2}$. Then $y$ is a $c$-median vertex since 
${\rm bw}_c(y,T) \leq N- c(y) \leq \frac{N}{2}$. It is
easy to see that $y$ is the only $c$-median vertex of $T$, so $y=v$. 
Now $v$ is an end-vertex of $T$. Indeed, if $v$ is an internal vertex, then $T-v$ has at least 
two components, and by Claim 1 their combined weight is more than $\frac{N}{2}$. Since 
$L>\frac{N}{2}$, this would imply that the total weight of $T$ is greater than $N$, a 
contradiction. Thus $v$ is an end-vertex.
By Claim 3 and Claim 4 it follows that $T$ 
is a path with all internal vertices having weight $k$.

We now show that $c(v)=L$. Indeed, if $c(v)>L$, then we obtain a weight function $c'$  from $c$ 
by reducing the weight of $v$ to $L$ and adding the excess weight $c(v)-L$ to $z_1$.  Then $c'$  
satisfies the hypothesis of the lemma and $v$ is a $c'$-median vertex of $T$. 
As above, $c'$ satisfies the hypothesis of the lemma, and we have
\[ \sigma_{c'}(v,T) = \sigma_{c}(v,T)  + d(v,z_1)(c(v)-L) > \sigma_{c}(v,T), \] 
a contradiction to the maximality of $\sigma_c(v,T)$.

We now show that $c(z_1)=k$. Suppose to the contrary that $c(z_1)>k$. Since $N-L$ 
is an integer multiple of $k$ and $c(v)=L$, it follows that $c(z_1)$ is an 
integer multiple of $k$, which implies
that $c(z_1) \geq 2k$. Extending the path $T$ by adding a new vertex $z_1'$
adjacent to $z_1$, and moving $c(z_1)-k$ weight units from $z_1$ to $z_1'$ yields a tree 
$H$ with a weight function $c'$ that satisfies the hypothesis of the lemma. Then 
$v$ is a $c'$-median vertex of $H$ and $\sigma_{c'}(v,H)> \sigma_c(v,T)$, a contradiction.

We have shown that $T$ is a path, the $c$-median vertex $v$ has weight $L$ and is an end-vertex 
of $T$, and all other vertices have weight $k$. 
So there are $\frac{N-L}{k}$ vertices of weight $k$ at distance 
$1,2,\ldots, \frac{N-L}{k}$ from $v$. Hence, 
\[ \sigma_{c}(v,T) = k\big(1 + 2 + \cdots + \frac{n-L}{k}  \big) 
                 = \frac{(N-L)(N-L+k)}{2k}, \]
which proves \eqref{eq:statement-of-main-lemma} for the case $L>\frac{N}{2}$. \\[1mm] 
{\sc Claim 6:} If $L \leq \frac{N}{2}$, then $v$ is an internal vertex of $T$, $T$ is a path, 
and $y$ is an end-vertex of $T$. \\
We first prove that $v$ is an internal vertex of $T$. Suppose to the contrary that $v$ is
an end-vertex. Then $T-v$ has only one component, and ${\rm bw}_c(v,T) = N- c(v)$. 
Since ${\rm bw}_c(v,T) \leq \frac{N}{2}$ it follows that $c(v) \geq \frac{N}{2}$.  
If this inequality is strict, i.e., if $c(v) > \frac{N}{2}$, then 
transferring $c(v)-\frac{N}{2}$ 
weight units from $v$ to its neighbour $u_1$ yields a weight function
$c'$ which satisfies the hypothesis of the lemma and for which $v$ is a $c'$-median
vertex. But $\sigma_{c'}(v,T) > \sigma_c(v,T)$, 
a contradiction to the maximality of $\sigma_c(v,T)$. Hence $c(v) = \frac{N}{2}$. 
Now consider vertex $u_1$. Clearly, the component of $T-u_1$
containing only $v$ has the maximum weight among all components of $T-u_1$, so
${\rm bw}_c(u_1,T) = \frac{N}{2}$, hence $u_1$ is also a $c$-median vertex. Hence $T$ 
has an internal vertex that is a $c$-median vertex. This contradicts the choice of $v$ as a 
$c$-median vertex that is also internal, if possible. It follows that our initial assumption 
that $v$ is an end-vertex is false, and so
$v$ is an internal vertex of $T$.
  
We now show that $T$ is a path. 
Suppose not. Since by Claim 3 each component of $T-v$ is a path where $v$ is adjacent to one of 
its ends, it follows that
$T-v$ has at least three components. Since $T-v$ has at most three components by Claim 2, 
it follows that $T-v$ has exactly three components.   

There are at least two components of $T-v$, $T_{i'}$ and $T_{i''}$ say, that do not contain $y$. 
Denote the third component by $T_i$. Recall that $z_{i'}$ and $z_{i''}$ are the end-vertices of 
$T$ in $T_{i'}$ and $T_{i''}$, respectively. We may assume that $d_T(v,z_{i'}) \geq d_T(v,z_{i''})$.            

Note that $c(z_{i'}) < 2k$ (and similarly, $c(z_{i''}) < 2k$) since otherwise, if  
$c(z_{i'})\geq 2k$, then we obtain a new graph by adding a new vertex of weight $k$, joining 
it to $z_{i'}$ and reducing the weight of $z_{i'}$ by $k$, which increases the weighted distance 
of $v$, which contradicts our choice of $T$ and $c$.

We now bound $C_{i'}$. Clearly, $C_{i'} = N - C_i - C_{i''}-c(v)$. Since 
$C_{i} + C_{i''} > \frac{N}{2}$ by Claim 1, and $c(v) \geq k$, it follows that 
$C_{i'} < \frac{N}{2} -k$. 

We may assume that $c(z_{i''})=k$. Indeed, if $c(z_{i''}) > k$, then let $r=c(z_{i''})-k$. As
shown above, we have $r \leq k$. Consider the weight function obtained from $c$ by shifting 
the extra weight $r$ from $z_{i''}$ to $z_{i'}$. This does not increase the total weight of 
$T_{i'}$ beyond $\frac{N}{2}$, so $v$ is also a $c'$-median 
vertex of $T$, and the weighted
distance of $v$ has not decreased since $d_T(v,z_{i'}) \geq d_T(v,z_{i''})$. Hence, from now 
on we may assume that $c(z_{i''})=k$.

Let $x$ be the neighbour of $z_{i''}$ in $T$. Consider the tree $H=T-z_{i''}x + z_{i''} z_{i'}$, 
so vertex $z_{i''}$ is transferred from $T_{i''}$ to $T_{i'}$. Since $C_{i'} < \frac{N}{2}-k$, 
vertex $v$ has $c$-branch weight at most $\frac{N}{2}$ in $H$ and is thus a $c$-median vertex of $H$. 
Moreover,
\[ \sigma_c(v,H) = \sigma_c(v,T) + k(d_T(v,z_{i'}) + 1 - d_T(v,z_{i''})) > \sigma_c(v,T), \]
a contradiction to the choice of $T$. Hence $T$ is a path. 

Now we complete the proof of Claim 6 by showing that $y$ is an end-vertex of $T$. It
suffices to show that $y \neq v$ since then by Claim 4 all internal vertices of $T$ have weight
$k$, so $y$ is an end-vertex.

Suppose to the contrary that $y=v$. Since $v$ is an internal vertex, it follows from
Claim 4 that $c(v)=L$. Let $c'$ be the weight function obtained from $c$ by moving $L-k$
weight units from $v$ to $u_2$, i.e., let $c'(v)=k$, $c'(u_2) = c(u_2) + L-k$ and 
$c'(x) = c(x)$ for all $x \in V(T)-\{ v,u_2\}$. Clearly, $c'$ satisfies the hypothesis of the lemma. 
We have either $C_1 > \frac{N}{2}-k$ or $C_1 \leq \frac{N}{2}-k$. 

If $C_1 > \frac{N}{2}-k$, then $v$ is also a $c'$-median vertex of $T$. Indeed, 
the weights with respect to $c'$ of the two branches of $T-v$ are $C_1$ and $N-C_1-k < \frac{N}{2}$.
Clearly, $\sigma_{c'}(v,T) = \sigma_c(v,T)+L-k > \sigma_c(v,T)$, a contradiction to the 
maximality of $\sigma_c(v,T)$.

If $C_1 \leq \frac{N}{2}-k$, then $u_2$ is a $c'$-median vertex of $T$. Indeed, 
the total weights with respect to $c'$ of the two branches of $T-u_2$ are $C_1+k$ and 
$C_2 - c(u_2)$, and both terms are clearly not more than $\frac{N}{2}$. 
Clearly, $\sigma_{c'}(u_2,T) =  \sigma_{c}(v,T) + C_1 + k - C_2 > \sigma_c(v,T)$, 
again a contradiction to the maximality of $\sigma_c(v,T)$. Hence Claim 6 holds. \\[1mm]
{\sc Claim 7:} If $L \leq \frac{N}{2}$, then the lemma holds. \\
Assume that $L \leq \frac{N}{2}$. Then by Claim 6, we have that $T$ is a path, $v$ is an internal
vertex, and $y$ is one of the two end-vertices of $T$. Denote the other end-vertex of $T$
by $x$. Moreover, denote the components of $T-v$ that contain vertex $x$ and $y$ by $T_x$ and
$T_y$, respectively. Note that $c(w)=k$ for all $w \in V(T)-\{x,y\}$ by Claim 4.
Let $C_x$ and $C_y$ be the total weight of the components of $T-v$ that contain vertex $x$ and
$y$, respectively. Let $r_x=c(x)-k$ and $r_y=c(y)-L$. We prove that
\begin{equation}   \label{eq:rx+ry} 
0 \leq r_x, r_y <k \quad \textrm{and} \quad r_x + r_y \in \{0,k\}. 
\end{equation} 
We first show that $0 \leq r_x < k$. Clearly, $0 \leq r_x$. 
If $r_x \geq k$, then we obtain a contradiction by adding a new vertex $x'$ and joining
it to $x$ to obtain a new tree $H$, and defining a new weight function $c'$ with $c'(x)=k$, 
$c'(x')=r_x$ and the remaining vertices have the same weight as for $c$. Then 
$v$ is a $c'$-median vertex of $H$ but $\sigma_{c'}(v,H) > \sigma_c(v,T)$, a contradiction. 
Similarly we show that $0 \leq r_y < k$. 
Since $N-L$ is an integer multiple of $k$, and since all vertices except possibly $x$ and $y$ have
weight $k$, it follows that $r_x+r_y$ is also an integer multiple of $k$. From $0 \leq r_x, r_y <k$
we conclude that $r_x+r_y$ either equals $0$ or $k$. This proves \eqref{eq:rx+ry}. 

We now consider three cases, depending on the values of $L$ and ${\rm bw}_c(v,T)$. \\[1mm]
{\sc Case 1:} $c(y)=L$. \\
Since $c(y)=L$ we have $r_y=0$ and thus $c(x)=k$ by \eqref{eq:rx+ry}.  Hence $T$ is a path 
$v_0, v_1,\ldots,v_{(N-L)/k}$, where $v_0$ has weight $L$ and the other vertices have weight $k$. 

In order to determine a $c$-median vertex of $T$, 
let $a \in \{1,2,\ldots, \frac{N-L}{k}\}$ be the largest value for which 
the component of $T-v_a$ containing $v_0$ has a total weight of
not more than $\frac{N}{2}$, so $a = \lfloor  \frac{N/2 - L + k}{k} \rfloor$. 
Then ${\rm bw}_c(v_a,T) \leq \frac{N}{2}$, and so $v_a$ is a 
$c$-median vertex of $T$. 
Apart from $v_a$ and $y$, $T$ has $a-1$ vertices of weight $k$ at distance
$1,2,\ldots,a-1$ from $v_a$ in the component of $T-v_a$ containing $v_0$, and 
$\frac{N-L - ak}{k}$ vertices of weight $k$ at distance $1,2,\ldots, \frac{N-L - ak}{k}$ 
from $v_a$ in the other component of $T-v_a$. Thus, 
\begin{eqnarray*} 
\sigma_{c}(v_a,T) & = & aL + k\big( 1 + 2+ \ldots + (a-1)\big) 
                        + k \Big( 1 + 2 + \ldots + \frac{N-L-ak}{k} \Big)  \\
          & = & aL + \frac{ka(a-1)}{2} + \frac{(N-L-ak)(N-L-ak+k)}{2k}.
\end{eqnarray*}     
Now $\frac{N-2L}{2k} < a \leq \frac{N-2L+2k}{2k}$ and so 
$N-L-ak < \frac{N}{2}$. Hence 
\begin{eqnarray*} 
\sigma_{c}(v_a,T) & < & \frac{N-2L+2k}{2k}L + \frac{k}{2} \frac{N-2L+2k}{2k}\frac{N-2L}{2k}
                        + \frac{1}{2k} \frac{N}{2} \Big( \frac{N}{2}+k\Big)  \\
          & = & \frac{N^2-2L^2}{4k} + \frac{N+L}{2},
\end{eqnarray*}  
and \eqref{eq:statement-of-main-lemma} holds in this case. \\[1mm]
{\sc Case 2:}  $c(y) > L$ and ${\rm bw}_c(v,T) = \frac{N}{2}$. \\
Since ${\rm bw}_c(v,T) = \frac{N}{2}$, there exists a component of $T-v$ with total
weight $\frac{N}{2}$. Let $v'$ be the neighbour of $v$ in this component. Then $T-vv'$ 
has two components, each of total weight $\frac{N}{2}$. It follows that 
${\rm bw}_c(v',T) = \frac{N}{2}$, and so $v'$ is also a $c$-median vertex. 
We assume that $x$ is in the same component of $T-vv'$ as $v$, and $y$ is in the same 
component as $v'$ (otherwise the proof is analogous). 
The weights of the two components of $T-vv'$ are $k d_T(v,x) + c(x)$ and $k d_T(v',y) + c(y)$. 
Since these weights are equal, and since $c(y) \geq c(x)$, we have
$k d_T(v,x) \geq k d_T(v',y)$ and thus $d_T(v,x) \geq d_T(v',y)$. Since $v$ and $v'$ are adjacent,
this implies that
\begin{equation} \label{eq:d(v,y) <= d(v,x)+1}
d_T(v,y) \leq d_T(v,x)+1.
\end{equation}
Now we obtain a new tree $H$ by adding a new vertex $x'$ and joining it to $x$. We define
a new weight function $c'$ by letting $c'(y)=L$, $c'(x)=k$, $c'(x')=r_x+r_y$, and as in 
$c$, all other vertices have weight $k$. By \eqref{eq:rx+ry} and $r_y>0$ we
have that $c'(x')=k$. Note that $v$ is a $c'$-median vertex. Indeed,  the 
components of $T-v$ containing $y$ and $x$ have (with respect to $c$) total weight $\frac{N}{2}$ 
and $\frac{N}{2}-k$, 
respectively, so the components of $H-x$ containing $y$ and $x$ have (with respect to $c'$) total weight 
$\frac{N}{2}-r_y$ and $\frac{N}{2}-k+r_y$, respectively, and both are less than $\frac{N}{2}$. 
From \eqref{eq:d(v,y) <= d(v,x)+1} we obtain that 
\[ \sigma_{c'}(v,H) = \sigma_c(v,T) + r_x + \big(d_T(v,x)+1 - d_T(v,y)\big)r_y > \sigma_c(v,T), \]
contradicting our choice of $v$, $T$ and $c$. Hence Case 2 cannot occur.  \\[1mm]
{\sc Case 3:}  $c(y) > L$ and ${\rm bw}_c(v,T) < \frac{N}{2}$. \\
We first show that $d_T(v,x)=d_T(v,y)$. Suppose not. Then we have either $d_T(v,x)<d_T(v,y)$ or 
$d_T(v,x)>d_T(v,y)$. First assume that $d_T(v,x) < d_T(v,y)$.  
Define $\varepsilon := \min\{\frac{N}{2} - {\rm bw}_c(v,T), c(y)-L, c(x)-k\}$. 
Since ${\rm bw}_c(v,T) < \frac{N}{2}$, we have that $\varepsilon >0$. 
We define a new weight function $c'$ by moving $\varepsilon$ weight units from $x$ to $y$. Clearly, 
$c'$ satisfies the hypothesis of the lemma and $v$ is a $c'$-median vertex of $T$. 
Then $\sigma_{c'}(v,T) > \sigma_c(v,T)$, a contradiction. We obtain a similar contradiction if $d_T(v,x)>d_T(v,y)$. 
This proves that $d_T(v,x)=d_T(v,y)$.

If $C_x \leq \frac{N}{2}-k$, then we obtain a new tree $H$ by adding a new vertex
$x'$ and joining it to $x$, and a new weight function $c'$ with $c'(x')=r_x+r_y=k$, $c(y)=L$ and
all remaining vertices have weight $k$.  Clearly, $v$ is a $c'$-median vertex since its 
$c'$-branch weight does not exceed $\frac{N}{2}$. Also, $\sigma_{c'}(v,H) > \sigma_c(v,T)$,
a contradiction. 
If $\frac{N}{2}-k < C_x < \frac{N}{2}$, then either $C_x \leq \frac{N}{2}-r_y$ or
$C_x > \frac{N}{2}-r_y$. If $C_x \leq \frac{N}{2}-r_y$, then we consider the weight function 
$c'$ obtained 
from $c$ by decreasing the weight of $y$ by $r_y$ and increasing the weight of $x$ by $r_y$.
Now $d_T(v,x) = d_T(v,y)$ implies that $\sigma_{c'}(v,T)=\sigma_c(v,T)$.  
Since $c'(y)=L$, we can now apply Case 1 to $c'$. 
If $C_x > \frac{N}{2}-r_y$, then  let $\varepsilon := \frac{N}{2}- C_x$. 
Consider the weight function $c'$ obtained from $c$ by decreasing the weight of $y$ by 
$\varepsilon$ and increasing the weight of $x$ by $\varepsilon$.
As above we have $\sigma_{c'}(v,T)=\sigma_c(v,T)$.  
Since ${\rm bw}_{c'}(v,T)= \frac{N}{2}$, we can now apply Case 2 to $c'$. 
In both cases, Claim 7 follows. This completes the proof of Lemma 
\ref{la:proximity-of-weighted}. \hfill $\Box$  \\

\begin{lemma}  \label{la:remoteness-of-weighted}
Let $G$ be a connected graph and $k,L \in \reals$ with $0 < k < L$. 
Let $c: V(G) \longrightarrow \reals^{\geq 0}$ be
a weight function with total weight $N$. Assume that $c(u) \geq k$ for 
every vertex $u \in V(G)$, and that $G$ contains a vertex of weight at least $L$. 
Assume further that $N-L$ is an integer multiple of $k$. If $v$ is a vertex of $G$, then
\begin{equation} \label{eq:statement-of-main-lemma-remoteness}
\sigma_c(v,G) \leq \frac{(N-L)(N+L-k)}{2k}.
\end{equation}        
\end{lemma}

{\bf Proof:}
We give only an outline the proof of Lemma \ref{la:remoteness-of-weighted} since it uses arguments
very similar to those in the proof of Lemma \ref{la:proximity-of-weighted}, but is much
less elaborate. 

Let $N,k,L$ be as in the hypothesis of the lemma. Let $G$, $v$ and $c$ be such that 
$\sigma_c(v,G)$ is maximised among all graphs, vertices and weight functions satisfying
the hypothesis of the lemma. We may assume that $G$ is a tree.

Then $v$ is an end-vertex of $G$ since otherwise, if $u_1$ and $u_2$ are two neighbours
of $v$, the graph $H=G-vu_2+u_1u_2$ satisfies $\sigma_c(v,H) > \sigma_c(v,G)$, a contradiction.
A similar argument shows that $G$ is a path. As in the proof of Lemma \ref{la:proximity-of-weighted}
we show that $v$ and all internal vertices of the path $G$ have weight equal to $k$, and that
the other end-vertex of $G$ has weight $L$. Thus $G$ has $\frac{N-L}{k}-1$ vertices at 
distance $1,2,\ldots,\frac{N-L-k}{k}$ from $v$, all of weight $k$, and one vertex
at distance $\frac{N-L}{k}$ from $v$ which has weight $L$. Hence
\[ \sigma_c(v,G) = k\Big(1+2+\cdots +\frac{N-L-k}{k}\Big) + \frac{N-L}{k}L 
                   = \frac{(N-L)(N+L-k)}{2k}, \]
which proves Lemma \ref{la:remoteness-of-weighted}. \hfill $\Box$

\section{Proximity and remoteness of graphs of given minimum degree and maximum degree} 
\label{section:Proximity and remoteness of graphs of given minimum degree and maximum degree}

In this section we present our main results. We give a bound on the proximity of
a graph in terms of order, minimum degree and maximum degree, and we construct graphs to 
show that our bound is sharp apart from an additive constant.

If $v$ is a vertex of $G$, then by $S(v)$ we mean the subgraph with vertex set 
$N_G(v)\cup \{v\}$ in which $v$ is adjacent to
each vertex of $N_G(v)$, and no other edges are present. 
We define the distance between a vertex $v$ and a set $B \subseteq V(G)$ as 
$\min_{u \in B} d(v,u)$. 

\begin{theorem}   \label{theo:proximity-in-terms-of-mindegree-maxdegree}
Let $G$ be a connected graph of order $n$, minimum degree $\delta$ and maximum degree $\Delta$. \\
If $\Delta > \frac{n}{2}-1$, then 
\[ \pi(G) \leq \frac{3(n-\Delta)^2}{2(n-1)(\delta+1)} + \frac{13}{2}. \]
If $\Delta \leq \frac{n}{2}-1$, then 
\[ \pi(G) \leq \frac{3(n^2 - 2\Delta^2)}{4(n-1)(\delta+1)} + \frac{35}{4}. \]
\end{theorem}

{\bf Proof:}
We first construct a spanning tree $T$ of $G$ as follows. Let $b_0$ be
a vertex of $G$ of maximum degree. Let $B_0 = \{b_0\}$ and $T_0 = S(b_0)$. Let $b_1$ be a
vertex at distance exactly $3$ from $B_0$, if one exists. Then there exists an edge $e_1$ joining
some vertex of $T_0$ to some vertex of $S(b_1)$. Let $T_1$ be the tree obtained from 
$T_0 \cup S(b_1)$ by adding the edge $e_1$ and let $B_1 = B_0 \cup \{b_1\}$. 
Let $b_2$ be a vertex at distance exactly $3$ from $B_1$, if one exists. Then there exists an 
edge $e_2$ joining some vertex of $T_1$ to some vertex of $S(b_2)$. Let $T_2$ be the tree 
obtained from $T_1 \cup S(b_2)$ by adding the edge $e_2$  and let $B_2 = B_1 \cup \{b_2\}$. 
Generally, for given $B_j$ and $T_j$, choose a vertex $b_{j+1}$ at
distance exactly $3$ from $B_j$, if one exists, let $e_j$ be an edge joining some vertex of $T_j$ 
to some vertex of $S(b_{j+1})$, let $T_{j+1}$ be the tree obtained from $T_j \cup S(b_{j+1})$ 
by adding the edge $e_{j+1}$ and let $B_{j+1} = B_j \cup \{b_{j+1}\}$. Repeat this procedure, 
say, for $r$ steps, until each vertex of $G$ is at distance at most $2$ from some 
vertex in $B_r$.

Let $T'=T_r$ and $B =  B_r =\{b_0, b_1,\ldots,b_r\}$. Then all vertices of $G$ are within 
distance at most $2$ from $B$ in $G$, and thus adjacent to 
some vertex in $T'$. Now we obtain a spanning tree $T$ of $G$ by joining every vertex that is not
in $T'$ to a neighbour in $T'$, which results in a spanning tree $T$ of $G$. 
Note that $T$ has the same maximum degree as $G$ 
since ${\rm deg}_G(b_0) = {\rm deg}_T(b_0)$.

Consider $T$ as a weighted tree in which every vertex of $T$ has weight $1$. Define
a new weight function on the vertices of $T$ by moving each weight to the nearest vertex 
in $B$. More
precisely, for every vertex $v$ of $T$ let $v_B$ be a vertex in $B$ closest to $v$ in $T$. 
We now move the weight of $v$ to $v_B$, that is, we define the weight function 
$c: V(T) \longrightarrow \reals^{\geq 0}$ by
\[ c(w) = | \{v \in V(T) \ | \ v_B=w\} | \]
for $w\in V(T)$. Every vertex of $T$ is within distance $2$ of some vertex in $B$,
hence we have $d(v,v_B) \leq 2$ for all $v\in V(T)-B$ and $d(v,v_B)=0$ for all $v\in B$. Hence, 
\begin{eqnarray*}
\sigma(w,T) & = & \sum_{v \in V(T)} d_T(w,v)  \\
            & \leq & \sum_{v \in V(T)} (d_T(w,v_B) + d_T(v,v_B) ) \\
            & \leq & \sum_{v \in V(T)} c(v) d_T(w,v) + 2(n-|B|) \\
            & \leq & \sum_{v \in V(T)} c(v) d_T(w,v) + 2(n-1).
\end{eqnarray*}
Thus, for all $w \in V(T)$, 
\begin{equation} \label{eq:sigma(w)-vs-sigma_c(w)}
\sigma(w,T) \leq \sigma_c(w,T) + 2(n-1).
\end{equation}
Since the weight of the vertices of $T$ is concentrated entirely in $B$, we have
\[ \sigma_c(w,T) = \sum_{v \in B} c(v) d(w,v). \]
Denote by $T^3$ the graph with vertex set $V(T)$ in which two vertices $u,v$ are adjacent if and only if $d_T(u,v)\leq 3$. Let $F$ be the subgraph of $T^3$ induced by $B$. 
It follows from the construction of $T$ and $B$ that 
for each vertex $b_i \in B$  with $i\geq 1$ there exists some vertex $b_j \in B$ with 
$j<i$ such that $d_T(b_i,b_j)=3$. This implies that for every $i \geq 1$  there exists 
a $(b_i, b_0)$-path in $F$, and so $F$ is connected.

We now consider $F$ with the weight function $c$. (Note that we use the notation $c$ also
for the restriction of $c$ to $B$.) Fix a $c$-median vertex $w_0$ of $F$. 
Since $d_T(w_0,v) \leq 3 d_F(w_0,v)$ for all $v \in B$, we have
\begin{equation}  \label{eq:T-vs-F} 
\sigma_c(w_0,T) \leq \sum_{v \in B}  3d_F(w_0,v) = 3\sigma_c(w_0,F). 
\end{equation}
Our aim is to bound the weighted distance
$\sigma_c(w_0,F)$ using Lemma \ref{la:proximity-of-weighted}. 
If $w \in B$, then every vertex $v \in N_G(w) \cup \{w\}$ satisfies $v_B=w$, hence we have
$c(w) \geq {\rm deg}_G(w)+1$. This implies that every vertex of $B$ has weight at least
$\delta+1$, and that $F$ contains a vertex, viz.\ $w_0$, of weight at least $\Delta+1$. 
Also, $c(B)=c(T)=n$.
However, $n-(\Delta+1)$ is not necessarily a multiple of $\delta+1$, so we apply 
Lemma \ref{la:proximity-of-weighted} not to $c$ but to a slightly modified weight function $c'$. 

There exists $q \in \reals$ with $0 \leq q  \leq \delta$ such that $n-(\Delta+1)+q$ 
is a multiple of $\delta+1$. 
Let $c'$ be the weight function obtained from $c$ by increasing the weight of $w_0$
by $q$ and leaving the other weights unchanged. Then $c'(F)=n+q$, each 
vertex of $F$ has weight at least $\delta+1$, and $F$ contains a vertex of weight at least
$\Delta+1$. Furthermore, $c'(F)-(\Delta+1)$ is a multiple of $\delta+1$. 
Clearly, $\sigma_c(w_0,T) \leq \sigma_{c'}(w_0,T)$ and $w_0$ is a $c'$-median vertex of $F$. 
For $N,L,k \in \reals^{\geq 0}$ let $f_1(N,L,k)  = \frac{(N-L)(N-L+k)}{2k}$ and 
$f_2(N,L,k) = \frac{N^2 - 2L^2}{4k} + \frac{N+L}{2}$. 
By Lemma \ref{la:proximity-of-weighted} we have 
\[ \sigma_{c}(w_0,F) \leq \sigma_{c'}(w_0,F) \leq \left\{ \begin{array}{cc} 
     f_1(n+q,\Delta+1, \delta+1) & \textrm{if $\Delta+1 > \frac{n+q}{2}$,} \\
     f_2(n+q,\Delta+1, \delta+1) & \textrm{if $\Delta+1 \leq \frac{n+q}{2}$.} 
         \end{array} \right. \]
In order to eliminate $q$, we observe that
if $\Delta+1 \leq \frac{n}{2}$, then $\Delta+1 \leq \frac{n+q}{2}$, and if 
$\Delta+1 \geq \frac{n+\delta+1}{2}$, then $\Delta+1 > \frac{n+q}{2}$. 
If $\frac{n}{2} < \Delta+1 < \frac{n+\delta+1}{2}$, then clearly 
$\sigma_{c'}(w_0,F) \leq \max\{ f_1(n+q,\Delta+1, \delta+1), f_2(n+q,\Delta+1, \delta+1) \}$.
It is easy to verify that $f_1(N,L,k) \geq f_2(N,L,k)$ if and only if 
$2L - 2\sqrt{kL} \leq N \leq 2L + 2\sqrt{kL}$. For $N=n+q$, $L=\Delta+1$ and 
$k=\delta+1$, this condition is satisfied if $\frac{n}{2} < \Delta+1 < \frac{n+\delta+1}{2}$. 
Hence we have $\sigma_{c'}(w_0) \leq f_2(n+q, \Delta+1, \delta+1)$ if 
$\frac{n}{2} < \Delta+1 < \frac{n+\delta+1}{2}$. It follows that 
\[ \sigma_{c}(w_0,F) \leq \left\{ \begin{array}{cc} 
     f_1(n+q,\Delta+1, \delta+1) & \textrm{if $\Delta+1 > \frac{n}{2}$,} \\
     f_2(n+q,\Delta+1, \delta+1) & \textrm{if $\Delta+1 \leq \frac{n}{2}$.} 
         \end{array} \right. \]
Since $f_1(N,L,k)$ and $f_2(N,L,k)$ are increasing in $N$, and since $n+q \leq n+\delta$, we 
obtain, after evaluating $f_1(n+\delta,\Delta+1, \delta+1)$ and $f_2(n+\delta,\Delta+1, \delta+1)$, 
\begin{equation} \label{eq:bound-on-sigma(w0)-both-cases}
\sigma_{c}(w_0,F)  \leq \left\{ \begin{array}{cc} 
     \frac{(n+\delta-\Delta-1)(n+2\delta-\Delta)}{2(\delta+1)}  
     & \textrm{if $\Delta+1 > \frac{n}{2}$,} \\
     \frac{(n+\delta)^2 - 2(\Delta+1)^2}{4(\delta+1)} + \frac{n+\delta +  \Delta+1}{2} 
     & \textrm{if $\Delta+1 \leq \frac{n}{2}$.} 
         \end{array} \right.  
\end{equation}         
{\sc Case 1:} $\Delta > \frac{n}{2}-1$. \\   
Consider the right hand side of \eqref{eq:bound-on-sigma(w0)-both-cases}. We have       
$(n+\delta-\Delta-1)(n+2\delta-\Delta) 
  = (n-\Delta)^2 + (n-\Delta)(3\delta-1) + 2\delta(\delta-1)$. Bounding 
$3\delta-1$ by $3(\delta+1)$ and $2\delta(\delta-1)$ by $3\delta(\delta+1)$ we obtain that 
$(n-\Delta)(3\delta-1) + 2\delta(\delta-1) \leq 3(\delta+1)(n-\Delta+\delta) \leq 3(\delta+1)(n-1)$,
and so
\begin{equation} \label{eq:bound-on-sigma(w0)-case1}
\sigma_{c}(w_0,F) \leq  \frac{(n-\Delta)^2}{2(\delta+1)} + \frac{3}{2}(n-1).
\end{equation}
Combining \eqref{eq:sigma(w)-vs-sigma_c(w)} for $w=w_0$, \eqref{eq:T-vs-F} and 
\eqref{eq:bound-on-sigma(w0)-case1} we obtain 
\begin{eqnarray*}
\sigma(w_0,T) & \leq & \sigma_{c}(w_0,T) + 2(n-1) \\
              & \leq &  3\sigma_{c}(w_0,F) + 2(n-1)  \\
              & \leq & \frac{3(n-\Delta)^2}{2(\delta+1)} + \frac{13}{2}(n-1).
\end{eqnarray*}
Division by $n-1$ now yields the bound in the theorem. \\[1mm]
{\sc Case 2:} $\Delta \leq \frac{n}{2}-1$. \\   
Consider the right hand side of \eqref{eq:bound-on-sigma(w0)-both-cases}. We have   
$(n+\delta)^2 - 2(\Delta+1)^2 = n^2-2\Delta^2 +2\delta n + \delta^2 - 4\Delta -2$. Bounding 
$2\delta n$ by $2(\delta+1)(n-1)$ and $\delta^2 - 4\Delta -2$ by $(\delta+1)(\delta-1)$, 
we obtain that 
$2\delta n + \delta^2 - 4\Delta -2 \leq (\delta+1)(2(n-1)+\delta-1)$. 
We thus obtain
\begin{eqnarray} 
\sigma_{c}(w_0,F) & \leq &  \frac{n^2 - 2\Delta^2}{4(\delta+1)} 
          + \frac{2n +  \Delta+\frac{3}{2}\delta- \frac{1}{2}}{2} 
                   \nonumber \\
 &\leq & \frac{n^2 - 2\Delta^2}{4(\delta+1)} + \frac{9}{4}(n-1), \label{eq:bound-on-sigma(w0)-case2}
\end{eqnarray}         
where in the last step we bounded $\Delta \leq n-1$ and $\delta \leq n-2$. 
Combining \eqref{eq:sigma(w)-vs-sigma_c(w)} for $w=w_0$, \eqref{eq:T-vs-F} and 
\eqref{eq:bound-on-sigma(w0)-case2} we obtain 
\begin{eqnarray*}
\sigma(w_0,T) & \leq & \sigma_{c}(w_0,T) + 2(n-1) \\
              & \leq &  3\sigma_{c}(w_0,F) + 2(n-1)  \\
              & \leq & \frac{n^2 - 2\Delta^2}{4(\delta+1)} + \frac{35}{4}(n-1).
\end{eqnarray*}
Since $\pi(G) \leq \pi(T) \leq \frac{1}{n-1} \sigma(w_0,T)$,  division by $n-1$ now yields 
the desired bound. \hfill $\Box$  \\

\begin{theorem}   \label{theo:remoteness-in-terms-of-mindegree-maxdegree}
Let $G$ be a connected graph of order $n$, minimum degree $\delta$ and maximum degree $\Delta$. 
Then  there exists a spanning tree $T$ of $G$ with
\[ \rho(T) \leq  \frac{3(n^2 -\Delta^2)}{2(n-1)(\delta+1)} + 7. \]
\end{theorem}

{\bf Proof:}
Let $B$, $T$, $F$, $c$, $q$ and $c'$  be as in the proof of 
Theorem \ref{theo:proximity-in-terms-of-mindegree-maxdegree}. 
Let $u$ be a vertex of maximum average distance in $T$, i.e., 
$\overline{\sigma}(u) = \rho(T)$. By the construction of $T$ there exists 
a vertex $u_B \in B$ with $d_T(u,u_B) \leq 2$. Hence 
\begin{equation} \label{eq:remoteness-T-1}
\sigma(u,T)   \leq \sigma(u_B,T) + 2(n-1).
\end{equation}
By \eqref{eq:sigma(w)-vs-sigma_c(w)}, and as in 
\eqref{eq:T-vs-F} we have 
\begin{equation} \label{eq:remotness-T-2} 
\sigma(u_B,T) \leq \sigma_c(u_B,T) + 2(n-1) \leq 3 \sigma_c(u_B,F) + 2(n-1). 
\end{equation}
Using the same arguments as in the proof of 
Theorem \ref{theo:proximity-in-terms-of-mindegree-maxdegree}, we bound $\sigma_c(u,F)$ with help 
of the weight function $c'$, to which we apply Lemma \ref{la:remoteness-of-weighted} 
(with $N=n+q$, $L=\Delta+1$ and $k=\delta+1$) and use the fact that $q \leq \delta$ to obtain that 
\begin{eqnarray*}
\sigma_c(u_B,F) & \leq & \sigma_{c'}(u_B,F) \\
    & \leq & \frac{(n+q-\Delta-1)(n+q+\Delta-\delta)}{2(\delta+1)}  \\
    & \leq & \frac{(n+\delta-\Delta-1)(n+\Delta)}{2(\delta+1)}.   
\end{eqnarray*}    
Now 
$(n+\delta-\Delta-1)(n+\Delta) 
  = (n-\Delta)(n+\Delta) + (\delta-1)(n+\Delta) 
  < n^2-\Delta^2 + (\delta+1)(2n-2)$, 
and so 
\begin{equation}
\sigma_c(u_B,F) <  \frac{(n-\Delta)(n+\Delta)}{2(\delta+1)} + n-1.  \label{eq:remoteness-F} 
\end{equation}
Combining \eqref{eq:remoteness-T-1}, \eqref{eq:remotness-T-2} and \eqref{eq:remoteness-F} and 
dividing by $n-1$ yields 
\begin{eqnarray*}
\sigma(u,T) & \leq  & \sigma(u_B,T) + 2(n-1) \\
     & \leq & 3 \sigma_c(u_B,F) + 4(n-1) \\
     & \leq & \frac{3(n^2-\Delta^2)}{2(\delta+1)} + 7(n-1).
\end{eqnarray*}
Since $\rho(T) = \frac{1}{n-1} \sigma(u,T)$, division by $n-1$ now yields the theorem. \hfill $\Box$ \\

Since $\rho(G) \leq \rho(T)$ for every spanning tree of a connected graph $G$, we have
the following corollary. 

\begin{corollary}   \label{coro:remoteness-in-terms-of-mindegree-maxdegree}
If  $G$ is a connected graph of order $n$, minimum degree $\delta$ and maximum degree $\Delta$. 
then 
\[ \rho(G) \leq  \frac{3(n^2 -\Delta^2)}{2(n-1)(\delta+1)} + 7. \]
\end{corollary}

\section{A sharpness example}
\label{section:A sharpness example}

We now construct a graph that shows that for fixed $\delta \geq 3$, and any given values of 
$n$ and $\Delta$ with $\delta< \Delta < n$, 
there are graphs of order $n$, minimum degree $\delta$ and maximum degree $\Delta$, whose 
proximity and remoteness are
within a constant of the bounds in Theorems \ref{theo:proximity-in-terms-of-mindegree-maxdegree}
and \ref{theo:remoteness-in-terms-of-mindegree-maxdegree}, respectively.

If $G_1, G_2,\ldots,G_k$ are graphs, then we define the {\em sequential sum}  
$G_1+G_2+\cdots + G_n$ to be the graph obtained from the disjoint union of the graphs
$G_1, G_2,\ldots,G_k$ by joining every vertex of $G_i$ to every vertex of $G_{i+1}$ 
for $i=1,2,\ldots,k-1$. By $K_n$ we mean the complete graph on $n$ vertices. 

Let $n, \Delta, \delta \in \naturals$ be given with $3 \leq \delta < \Delta <n$. 
For the following construction we assume that $n-\Delta$ is a multiple of $\delta+1$, but it 
is not hard to modify the construction to work without this additional assumption. 
Let $k := \frac{n-\Delta}{\delta+1}$. Define the graph $G_{n,\Delta, \delta}$ by 
\[ G_{n,\Delta,\delta} = K_{\delta} + K_1 + [K_1 + K_{\delta-1} + K_1]^{k-1} + K_1 +K_{\Delta-1}, \]
where $[K_1 + K_{\delta-1} + K_1]^{k-1}$ stands for $k-1$ repetitions of the pattern 
$K_1 + K_{\delta-1} + K_1$. Clearly, $G_{n, \Delta, \delta}$ has order $n$, 
minimum degree $\delta$ and maximum degree $\Delta$.

We first bound the proximity of $G_{n,\Delta,\delta}$ from below. 
Define $V_1$ to be the set of $\delta+1$ vertices that belong to the first (counted from the left) or 
second complete graph in the sequential sum, i.e, to  $K_{\delta}\cup K_1$. 
For $i=1,2,\ldots,k-1$ let $V_{i+1}$ be the set of $\delta+1$ vertices that belong to $i$th 
repetition of the pattern $K_1 + K_{\delta-1} + K_1$ in the definition of $G_{n,\Delta,\delta}$.
Define $V_{k+1}$ to be the set of $\Delta$ vertices belonging to one of the rightmost two complete
graphs, i.e., to $K_1 +K_{\Delta-1}$ in the definition of $G_{n,\Delta,\delta}$. 
We make use of the fact that, whenever $x$ and $y$ are two vertices with $x\in V_i$ and 
$y \in V_j$, then $d(x,y) \geq 3|i-j|-2$. \\[1mm]
{\sc Case 1:}  $\Delta \leq \frac{n}{2}$. \\
Let $a:=\lceil \frac{n/2}{\delta+1} \rceil$. Then it is easy to see (for example using the fact 
that both components of $G_{n,\Delta,\delta} - V_a$ contain not more than $\frac{n}{2}$
vertices), that $V_a$ contains a median vertex, $w$ say. 
We bound $\sigma(w,G_{n,\Delta,\delta})$ from below by taking into account only the distances to vertices not in $V_a$.
\begin{eqnarray*}
\sigma(w,G_{n,\Delta,\delta}) & > & \sum_{i=1}^{a-1} \sum_{v \in V_i} d(w,v) + 
           \sum_{i=a+1}^{k} \sum_{v \in V_i} d(w,v) +  \sum_{v \in V_{k+1}} d(w,v) \\
    & \geq & \sum_{i=1}^{a-1} (\delta+1) [3(a-i)-2] + 
           \sum_{i=a+1}^{k} (\delta+1) [3(i-a)-2] \\
    & &  +  \Delta[3(k-a+1)-2] \\           
    & = & \frac{3}{2}(\delta+1) a (a-1) - 2(\delta+1)(a-1) + 
        \frac{3}{2}(\delta+1) (k-a+1) (k-a) \\
    & & - 2(\delta+1)(k-a) + 3\Delta (k-a+1) - 2\Delta \\
    & = & \frac{3}{2}(\delta+1) [a^2 + (k-a+1)^2] + 3\Delta(k-a+1) \\
    & &       - \big[(\delta+1)(\frac{7}{2}k-\frac{1}{2})+2\Delta\big].  
\end{eqnarray*}
Using the inequalities $a=\lceil \frac{n/2}{\delta+1} \rceil \geq \frac{n/2}{\delta+1}$, 
$k-a+1 = k - \lceil \frac{n/2}{\delta+1} \rceil + 1 > k - \frac{n/2}{\delta+1} 
     = \frac{n - 2\Delta}{2(\delta+1)}$, and 
     $(\delta+1)(\frac{7}{2}k-\frac{1}{2})+2\Delta < \frac{7}{2}(n-1)$ we obtain that
\begin{eqnarray*}
\sigma(w,G_{n,\Delta,\delta}) & > & \frac{3}{2} (\delta+1) 
         \big[ \frac{n^2}{4(\delta+1)^2} + \frac{(n-2\Delta)^2}{4(\delta+1)^2} \big]
         + 3\Delta \frac{n-2\Delta}{2(\delta+1)} - \frac{7}{2}(n-1) \\
         & = & \frac{3(n^2-2\Delta^2)}{4(\delta+1)} - \frac{7}{2}(n-1).
\end{eqnarray*}
Division by $n-1$ yields that
\[ \pi(G_{n,\Delta,\delta}) = \overline{\sigma}(w,G_{n,\Delta,\delta}) 
       > \frac{3(n^2-2\Delta^2)}{4(n-1)(\delta+1)} - \frac{7}{2}, \]
and so   $\pi(G_{n,\Delta,\delta})$ differs from the bound in  
Theorem \ref{theo:proximity-in-terms-of-mindegree-maxdegree} by less than $\frac{49}{4}$.
 \\[1mm]
{\sc Case 2:} $\Delta \geq \frac{n}{2}$. 
Clearly, $V_{k+1}$ contains a median vertex, $w$ say. We bound $\sigma(w, (G_{n,\Delta,\delta})$ 
from below by taking into account only the distances to vertices not in $V_{k+1}$.
\begin{eqnarray*}
\sigma(w, (G_{n,\Delta,\delta})) & > & \sum_{i=1}^{k} \sum_{v \in V_i} d(w,v)  \\
    & \geq & \sum_{i=1}^{k} (\delta+1) [3(k+1-i)-2] \\           
    & = & \frac{3}{2}(\delta+1) k (k+1) - 2(\delta+1)k  \\
    & = & \frac{3}{2}(\delta+1) k^2 - \frac{1}{2} (\delta+1)k.  
\end{eqnarray*}
Substituting $k= \frac{n-\Delta}{\delta+1}$, bounding 
$\frac{1}{2} (\delta+1)k = \frac{1}{2}(n-\Delta) \leq \frac{n-1}{2}$
and dividing by $n-1$ we obtain 
\[ \pi(G_{n,\Delta,\delta}) = \overline{\sigma}(w, (G_{n,\Delta,\delta}) 
       > \frac{3(n-\Delta)^2}{(n-1)(\delta+1)} - \frac{1}{2}, \]
and so $\pi(G_{n,\Delta,\delta})$  differs from the bound in  
Theorem \ref{theo:proximity-in-terms-of-mindegree-maxdegree} by less than 
$6 \delta+\frac{5}{2}$, which for fixed $\delta$ is a constant. \\

We now bound the remoteness of $G_{n,\Delta,\delta}$. Let $u$ be a vertex of the
graph $K_{\delta}$ in the representation of $G_{n,\Delta,\delta}$ as a sequential sum. 
 We have 
\begin{eqnarray*}
\sigma(u, (G_{n,\Delta,\delta})) & = & \sum_{v\in V_1} d(u,v) + \sum_{i=2}^{k} \sum_{v \in V_i} d(u,v) 
    + \sum_{v\in V_{k+1}} d(u,v)    \\
    & = & \delta + \sum_{i=2}^{k} (\delta+1) [3i-3] + \Delta\, 3k -1 \\           
    & = & 3(\delta+1) \frac{k(k-1)}{2} + 3\Delta k + \delta-1. 
\end{eqnarray*}
Substituting  $k = \frac{n-\Delta}{\delta+1}$ yields after simplification that
\[ \sigma(u, (G_{n,\Delta,\delta})) = \frac{3(n^2 -\Delta^2)}{2(\delta+1)} - \frac{3}{2}(n- \Delta - \frac{2}{3}(\delta+1)) 
            > \frac{3(n^2 -\Delta^2)}{2(\delta+1)} - \frac{3}{2}(n-1).
 \]
Division by $n-1$ yields that 
\[ \rho(G_{n,\Delta,\delta}) \geq \overline{\sigma}(u, (G_{n,\Delta,\delta}) 
    > \frac{3}{2} \frac{n^2 -\Delta^2}{(n-1)(\delta+1)} - \frac{3}{2}, \] 
and so $\rho(G_{n,\Delta,\delta})$ differs from the bound in 
Corollary \ref{coro:remoteness-in-terms-of-mindegree-maxdegree} by not more than $\frac{17}{2}$.



\end{document}